\documentclass[11pt,a4paper]{article}

\usepackage{theorem,enumerate}
\usepackage{amsmath,latexsym,amssymb,amsfonts}
\usepackage{eucal}
\usepackage{mathrsfs}
\usepackage{array}
 \usepackage{graphicx, color}

\newcounter{Scounter}
\newtheorem{theorem}{Theorem}
\newtheorem{lemma}[theorem]{Lemma}

\newtheorem{corollary}[theorem]{Corollary}

\newtheorem{claim}[theorem]{Claim}

\newtheorem{definition}[theorem]{Definition}

\newtheorem{remark}[theorem]{Remark}

\numberwithin{theorem}{section}

\newcommand{\Proof}{\noindent\textbf{Proof}.\quad}

\newcommand{\Case}[1]{\noin{\bf Case #1}}

\newcommand{\qed}{{$\quad\square$\vs{3.6}}}

\newcommand{\vs}[1]{\vspace*{#1 mm}}

\newcommand{\noin}{\noindent}

\def\M{{ \mathcal{M}}}

\bfseries\normalfont

\makeatletter
\def\thanks#1{%
   \footnotemark
   \edef\@tempa{\noexpand\noexpand\noexpand\footnotetext[\the\c@footnote]}%
   \toks@\expandafter{\@thanks}%
   \toks\tw@{{#1}}
   \xdef\@thanks{\the\toks@\@tempa\the\toks\tw@}}
\makeatother

\begin{document}

\title{Planar graphs without 4-cycles adjacent to triangles are DP-4-colorable}

\author{
Seog-Jin Kim\thanks{Department of Mathematics Education, Konkuk University,
Korea. Email: {\tt skim12@konkuk.ac.kr}\\
Supported by Basic Science Research Program through the National Research Foundation of Korea(NRF) funded by the Ministry of Education(NRF-2015R1D1A1A01057008).}
and
Xiaowei Yu\thanks{School of Mathematics, Shandong University, Jinan, Shandong, 250100, P. R. China. Corresponding author.
Email: {\tt xwyu2013@163.com}\\Supported by the National Natural Science Foundation of China (11371355,
11471193, 11271006, 11631014), the Foundation for Distinguished Young Scholars of Shandong Province
(JQ201501), and fundamental research funding of Shandong University.
}
}

\maketitle

\begin{abstract}
DP-coloring (also known as correspondence coloring) of a simple graph is
a generalization of list coloring.
It is known that planar graphs without 4-cycles adjacent to triangles are 4-choosable, and planar graphs without 4-cycles are DP-4-colorable. In this paper, we show that  planar graphs without 4-cycles adjacent to triangles are DP-4-colorable, which is an extension of the two results above.
\end{abstract}


\section{Introduction}


We use standard notation. For a set $S$, Pow$(S)$ denote the power set of $S$, i.e., the set of all subsets of $S$. We denote by $[k]$ the set of integers from $1$ to $k$. All graphs considered here are finite, undirected, and simple. For a graph $G$, $V(G)$, $E(G)$, and $F(G)$ denote the vertex sets, edge sets and face sets of $G$, respectively. For a set $U\subseteq V(G)$, $G[U]$ is the subgraph of $G$ induced by $U$.

Recall that a \emph{proper $k$-coloring} of a graph $G$
is a mapping $f : V(G) \rightarrow [k]$
such that $f(u) \not= f(v)$ for any $uv \in E(G)$.
The minimum integer $k$
such that $G$ admits a proper coloring
is called the \emph{chromatic number of $G$},
and denoted by $\chi(G)$.

List coloring is a generalization of graph coloring that was introduced independently by Vizing \cite{Vizing} and Erd\H{o}s, Rubin, and Taylor \cite{ERT}.
Let $C$ be a set of colors.
A \emph{list assignment} $L : V(G) \rightarrow Pow(C)$ of $G$
is a mapping that assigns a set of colors to each vertex. If $|L(v)|\ge k$ for all $v\in V(G)$, then $L$ is called a $k$-\emph{list assignment}.
A proper coloring $f: V(G) \rightarrow C$ is called an \emph{$L$-coloring} of $G$ if
$f(u) \in L(u)$ for any $u \in V(G)$. The \emph{list-chromatic number} or the \emph{choice number} of $G$, denoted by $\chi_{\ell}(G)$, is the smallest $k$ such that $G$ admits an $L$-coloring for every $k$-list assignment $L$ for $G$.


%
%

Since a proper $k$-coloring corresponds to an $L$-coloring
with $L(u) = [k]$ for any $u \in V(G)$,
we have $\chi(G) \leq \chi_{\ell}(G)$.
It is well-known that
there are infinitely many graphs $G$
satisfying $\chi (G) < \chi_{\ell}(G)$,
and the gap can be  arbitrarily large.


\bigskip

In order to consider some problems on list chromatic number,
Dvo\v{r}\'{a}k and Postle \cite{DP} considered
a generalization of a list coloring.
They call it a \emph{correspondence coloring},
but we call it a \emph{DP-coloring} for short,
following Bernshteyn, Kostochka and Pron \cite{BKP}.

Let $G$ be a  graph
and $L$ be a list assignment of $G$.
For each edge $uv$ in $G$,
let $M_{L,uv}$ be an arbitrary matching (maybe empty)
between $\{u\} \times L(u)$ and $\{v\} \times L(v)$.
Without abuse of notation,
we sometimes regard $M_{L,uv}$
as a bipartite graph in which the edges are
between $\{u\} \times L(u)$ and $\{v\} \times L(v)$, and
the maximum degree is at most 1.

\begin{definition}
Let $\M_{L} = \big\{M_{L,uv} : uv \in E(G) \big\}$,
which is called a \emph{matching assignment over $L$}.
Then a graph $H$ is said to be the \emph{$\M_{L}$-cover} of $G$
if it satisfies all the following conditions:
\begin{enumerate}[{\upshape (i)}]
\item
The vertex set of $H$ is
$\bigcup_{u \in V(G)} \big(\{u\} \times L(u)\big)
= \big\{(u,c): u \in V(G), \ c \in L(u)\big\}$.

\item
For every $u \in V(G)$,
the graph $H[\{u\} \times L(u)]$ is a clique.
\item
For any edge $uv$ in $G$,
$\{u\} \times L(u)$
and
$\{v\} \times L(v)$
induce the graph obtained from $M_{L,uv}$ in $H$.

\end{enumerate}
\end{definition}

\begin{definition}
An $\M_{L}$-coloring of $G$ is an independent set $I$ in
the $\M_{L}$-cover with $ |I| = |V(G)|$.
The \textit{DP-chromatic number},
denoted by $\chi_{\text{DP}}(G)$,
is the minimum integer $k$
such that $G$ admits an $\M_{L}$-coloring
for each $k$-list assignment $L$ and each matching assignment $\M_{L}$ over $L$.  We say that a graph $G$ is \emph{DP-$k$-colorable} if $\chi_{DP}(G) \leq k$.
\end{definition}

Note that
when $G$ is a simple graph and
$$M_{L,uv} = \big\{(u,c)(v,c): c \in L(u) \cap L(v)
\big\}
$$
for any edge $uv$ in $G$,
then $G$ admits an $L$-coloring
if and only if
$G$ admits an $\M_{L}$-coloring.
This implies $\chi_{\ell}(G) \leq \chi_{\text{DP}}(G)$. Thus DP-coloring is a generalization of the list coloring. Thus, given the fact $\chi_{\ell} (G) \leq k$, it is interesting to check whether $\chi_{DP}(G) \leq k$ or not.

Dvo\v{r}\'{a}k and Postle \cite{DP} showed that $\chi_{DP}(G) \leq 5$ if $G$ is a  planar graph, and $\chi_{DP}(G) \leq 3$ if $G$ is a planar graph with girth at least 5.
Also, Dvo\v{r}\'{a}k and Postle \cite{DP} observed that $\chi_{DP}(G) \leq k+1$ if $G$ is $k$-degenerate.

On the other hand, there are some differences between DP-coloring and list coloring.
There are infinitely many simple graphs $G$
satisfying $\chi_{\ell} (G) < \chi_{\text{DP}}(G)$:
It is known that
$\chi(C_n) = \chi_{\ell}(C_n) = 2 < 3 = \chi_{\text{DP}}(C_n)$
for each even integer $n \geq 4$ (see \cite{BKP}).
Furthermore,
the gap $\chi_{\text{DP}}(G) - \chi_{\ell} (G)$
can be arbitrary large.
For example,
Bernshteyn \cite{Bernshteyn} showed that
for a simple graph $G$ with average degree $d$,
we have $\chi_{\text{DP}}(G) = \Omega(d / \log d)$,
while Alon \cite{Alon} proved
that $\chi_{\ell}(G) = \Omega(\log d)$ and
the bound is sharp.
Recently, there are some other works on DP-colorings,
see \cite{ BK2, BKZ, KO}.


\bigskip

Thomassen \cite{Thomassen} showed that every planar graph is 5-choosable, and Voigt \cite{Voigt} showed that there are planar graphs which are not 4-choosable.
Thus finding sufficient conditions for planar graphs to be 4-choosable is an interesting problem.

Two faces are {\em adjacent} if they have at least one common edge, and two faces are {\em normally adjacent} if they are adjacent and have exactly one common edge.
Let $C_k$ be the cycle of length $k$.
Lam, Xu, and Liu \cite{Lam} verified that every planar graph without $C_4$ is 4-choosable.  And Cheng, and Chen, Wang \cite{CCW}, and Kim and Ozeki \cite{KO2} extended the result independently by certifying the following theorem.

\begin{theorem}\label{listversion} The following results hold independently.
\begin{description}
  \item[A.] (\cite{CCW}) If $G$ is a planar graph without $4$-cycles adjacent to $3$-cycles, then $\chi_{\ell}(G)\le 4$.
  \item[B.] (\cite{KO2}) If $G$ is a planar graph without $4$-cycles, then $\chi_{DP}(G)\le 4$.
\end{description}
\end{theorem}

In this paper, we extend Theorem \ref{listversion} by proving
the following theorem.

\begin{theorem} \label{main-thm} If $G$ is a planar graph without $4$-cycles adjacent to $3$-cycles, then $\chi_{DP} (G) \leq 4$.
\end{theorem}

Even though $\chi_{DP} (G)=\chi_{\ell}(G)$ for some special graphs, Theorem \ref{main-thm} is not trivial from Theorem \ref{listversion} (A). We will give an explanation in section 3.1.



\section{Proof of Theorem \ref{main-thm}}

Suppose that Theorem \ref{main-thm} does not hold.
In the rest of paper, let $G$ be a minimal counterexample to Theorem \ref{main-thm} with fewest edges. From our hypothesis, graph $G$ has the following properties:
\begin{enumerate}[(a)]
\item $G$ is connected; and
\item $G$ has no subgraph isomorphic to 4-cycles adjacent to 3-cycles; and
\item $G$ is not DP-$4$-colorable; and
\item any proper subgraph $G'$ of $G$ is DP-$4$-colorable.
\end{enumerate}

Embedding $G$ into a plane, we obtain a plane graph $G=(V,E,F)$
where $V,E,F$ are the sets of vertices, edges, and faces of $G$,
respectively.

For a vertex $v\in V$, the \emph{degree} of $v$ in $G$ is
denoted by $d_G(v)$.  A vertex of degree $d$ (at least $d$, at most $d$, respectively) is called a $d$-\textit{vertex} ($d^+$-\textit{vertex}, $d^{-}$-\textit{vertex}, respectively). The notions of $d$-\textit{face}, $d^+$-\textit{face}, $d^{-}$-\textit{face} are similarly defined. According to (b), if a $5$-face is adjacent to a $3$-face, then they are normally adjacent.

For a face $f\in F$, if the vertices on $f$ in a cyclic order are $v_1, v_2, \ldots, v_k$, then we write $f= [v_1 v_2 \cdots v_kv_1]$, and call $f$ a $(d_G(v_1), d_G(v_2), \ldots, d_G(v_k) )$-face.

\subsection{Structures}

Using the properties of $G$ above, we can obtain several local structures of $G$.

\begin{lemma} \label{basic-lemma} Graph $G$ has no $3^-$-vertex.
\end{lemma}

\Proof
Suppose to the contrary that there exists a $3^-$-vertex $w$ in $G$. Let $L$ be a list assignment of $G$ with $|L(v)|\ge 4$ for any $v\in V$, and let $\M_{L}$ be a matching assignment over $L$. Let $G': = G -\{ w\}$ and $L'(v)=L(v)$ for $v\in V(G')$. According to (d), $G'$ admits an $\M_{L'}$-coloring.
Thus there is an independent set $I'$ in $\M_{L'}$-cover with $|I'| = |V(G)| - 1$.
For $w$, we define that
\[
L^*(w) = L(w) \setminus \bigcup_{uw \in E(G)} \big\{c' \in L(w): (u,c)(w,c') \in M_{L,uw}  \mbox{ and } (u, c) \in I' \big\}.
\]
Since $|L(w)| \ge 4 $ and $w$ is a $3^-$-vertex, we have that $|L^*(w)| \geq 1$. We denote by $\M_{L^*}$ the restriction of $\M_{L}$ into $G[w]$ and $L^*$. Obviously, $I=I'\cup \{\{w\}\times L^*(w)\}$ is an independent set in the $\M_{L^*}$-cover with $|I| = |V(G)|$, which is a contradiction to (c).
\qed

\begin{figure}[htbp]
 \begin{center}
\scalebox{0.65}[0.65]{\includegraphics {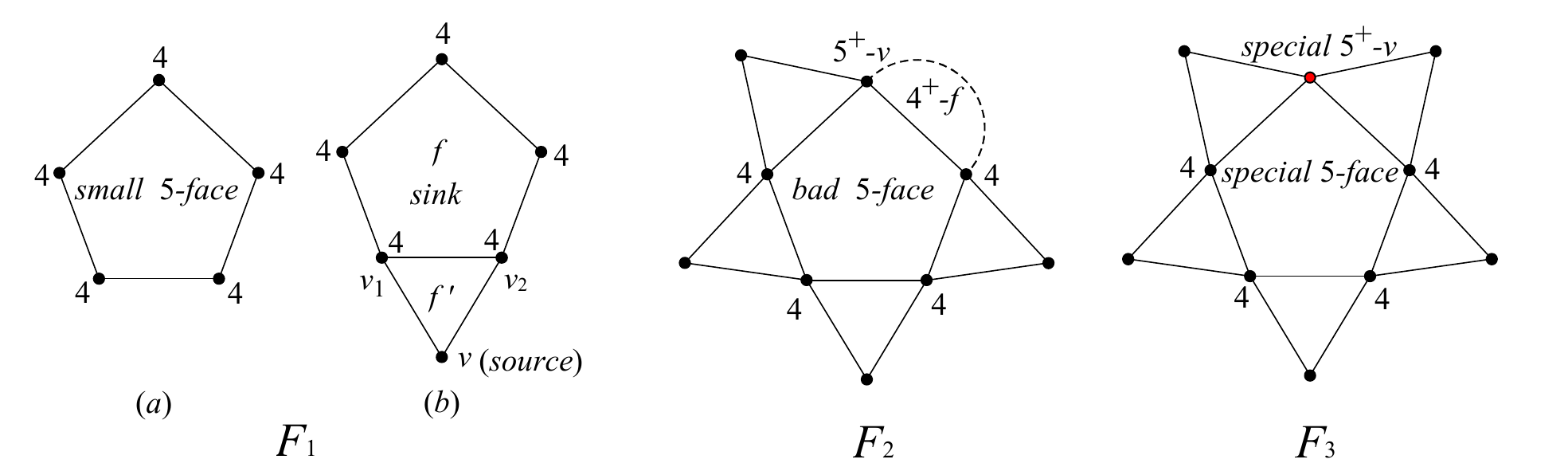}}
 \end{center}
 \caption{Definition of different 5-faces} \label{fig2-1}
\end{figure}

 Next we define several different $5$-faces.
Let $f$ be a $5$-face of $G$.

\begin{enumerate}[(i)]

\item If $f$ is a $(4,4,4,4,4)$-face, then we call $f$ a \emph{small} 5-face. Assume that $f$ is adjacent to a $3$-face $f'=[vv_1v_2v]$ with the common edge $v_1v_2$, then we call $v$ a \emph{source} of $f$. Equivalently, the face $f$ is called a \emph{sink} of $v$ (see $F_1$ in Figure \ref{fig2-1}).

\item If $f$ is a $(5^+,4,4,4,4)$-face, and $f$ is incident to four $3$-faces and one $4^+$-face, and the $4^+$-face is incident to the $5^+$-vertex on $f$, then we call $f$ a \emph{bad} $5$-face (see $F_2$ in Figure \ref{fig2-1}).

\item If $f$ is a $(5^+,4,4,4,4)$-face, and $f$ is incident to five $3$-faces, then we call $f$ a \emph{special} $5$-face. Meanwhile, we call the $5^+$-vertex on $f$ a \emph{special} vertex (see $F_3$ in Figure \ref{fig2-1}).
\end{enumerate}

\begin{remark}\label{remark2}
\begin{enumerate}[(1)]
\item A special vertex is incident to at least one special $5$-face, and a special $5$-face is incident to exactly one special vertex.
\item If there exist two $5^+$-vertices on a $5$-face $f$, then $f$ is neither special nor bad.
\end{enumerate}
\end{remark}

\begin{lemma} \label{good vertex}
Every source is a $5^{+}$-vertex.
\end{lemma}
\Proof
Let $f = [v_1 v_2 v_3 v_4 v_5 v_1]$ be a small 5-face, and let $z$ be a source of $f$. By Lemma \ref{basic-lemma}, there exists no $3^-$-vertex in $G$, thus we suppose that $d_G(z)  = 4$.
Without loss of generality, we may assume that $z$ is adjacent to $v_1$ and $v_2$.
Let $G[S]$ be the subgraph of $G$ induced by $S = \{z, v_1, v_2, v_3, v_4, v_5 \}$.
Let $L$ be a list assignment of $G$ with $|L(v)|\ge 4$ for all $v\in V$,
and let $\M_{L}$ be a matching assignment over $L$.
Consider subgraph $G' := G - V(G[S])$ and $L'(v) = L(v)$ for $v \in V(G')$.
By (d),
$G'$ admits an $\M_{L'}$-coloring.
Thus there is an independent set $I'$ in the $\M_{L'}$-cover with $|I'| = |V(G)| - |S| = |V(G)| - 6$.  For $v \in \{v_1, v_2, v_3, v_4, v_5, z \}$, we define
\[
L^*(v) = L(v) \setminus \bigcup_{uv \in E(G)} \big\{c' \in L(v): (u,c)(v,c') \in M_{L,uv}  \mbox{ and } (u, c) \in I' \big\}.
\]
Because $|L(v)| \ge 4 $ for all $v \in V(G)$,
we have that $|L^*(v_1)| \geq 3, \ |L^*(v_2)| \geq 3$, and $|L^*(v)| \geq 2$ for $v \in \{v_3, v_4, v_5, z \}$.
We denote by $\M_{L^*}$
the restriction of $\M_{L}$ into $G[S]$ and $L^*$.

\begin{claim} \label{key-claim}
The $\M_{L^*}$-cover has an independent set $I^*$ with $|I^*| = 6 = |V(G[S])|$.
\end{claim}
\Proof Because $d_G(v)=4$ for $v\in S$, it holds that $|L^*(v_1)| \geq 3$ and $|L^*(z)| \geq 2$, we can color $v_1$ by $c \in L^*(v_1)$ such that $L^*(z) \setminus \{c' : (v_1, c) (z, c') \in M_{L^*}\}$ has at least two available colors.
By coloring greedily $v_5, v_4, v_3, v_2, z$ in order, we can find an independent set $I^*$ with $|I^*| = 6$.
This completes the proof of Claim \ref{key-claim}. \qed

Thus $G$ admits an $\M_L$-coloring $I = I' \cup I^*$ such that $|I| = |I'| + |I^*| = |V(G)|$.  It implies that $G$ is DP-4-colorable.  This is contradiction to (c).  Thus $d_G(z) \geq 5$.\qed

\begin{lemma}\label{3-faces} Let $v$ be a $5^+$-vertex of $G$. Then the following hold:

(1) $v$ is incident to at most $\lfloor\frac{d(v)}{2}\rfloor$ $3$-faces; and

(2) if $v$ is incident to $t$ $3$-faces and $2t<d_G(v)$, then $v$ is incident to at most $t-1$ special $5$-faces.
\end{lemma}

\Proof
Note that (1) holds obviously from (b).

And from the definition of a special $5$-face, if $v$ is a $5^+$-vertex of $G$, then each of the special $5$-faces incident to $v$ is adjacent to two $3$-faces, and both of the two $3$-faces are incident to $v$, which implies (2).
\qed

\begin{lemma}\label{specialgood}
Let $v$ be a $5^+$-vertex of $G$. Assume that $v$ is incident to a special or bad $5$-face $f_1$ and a $3$-face $f_2$ such that $f_1$ and $f_2$ are adjacent, then $v$ has no sink adjacent to $f_2$.
\end{lemma}

\begin{figure}[htbp]
 \begin{center}
\scalebox{0.6}[0.6]{\includegraphics {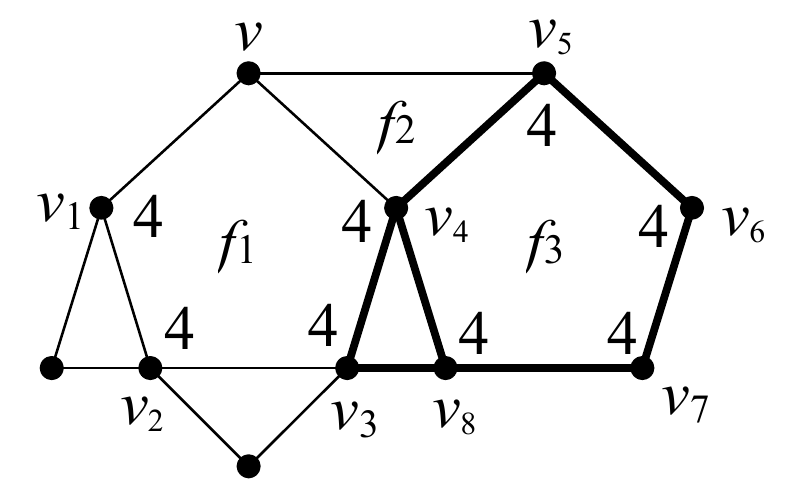}}
 \end{center}
 \caption{Illustration of Lemma \ref{specialgood}} \label{fig2-2}
\end{figure}

\Proof
Let $f_1=[vv_4v_3v_2v_1v]$ be a special or a bad $5$-face. According to the definition of special $5$-face and bad $5$-face, there are at least four 3-faces adjacent to $f_1$. Assume that $f_2=[vv_5v_4v]$, and $f_2$ is adjacent to $f_1$ (see Figure \ref{fig2-2}). Suppose to the contrary that $v$ has a sink $f_3=[v_4v_5v_6v_7v_8v_4]$ adjacent to $f_2$ with a common edge $v_4v_5$, then $v_3$ is a source of $f_3$. However, $v_3$ is a $4$-vertex, contrary to Lemma \ref{good vertex}.
\qed
\begin{lemma}\label{bad5face1} Let $f_1$ and $f_2$ be two bad $5$-faces. Then they can not normally adjacent with one common edge $vv_1$, where $v$ is the $5^+$-vertex on $f_1$ and $f_2$.

\end{lemma}
\begin{figure}[htbp]
 \begin{center}
\scalebox{0.6}[0.6]{\includegraphics {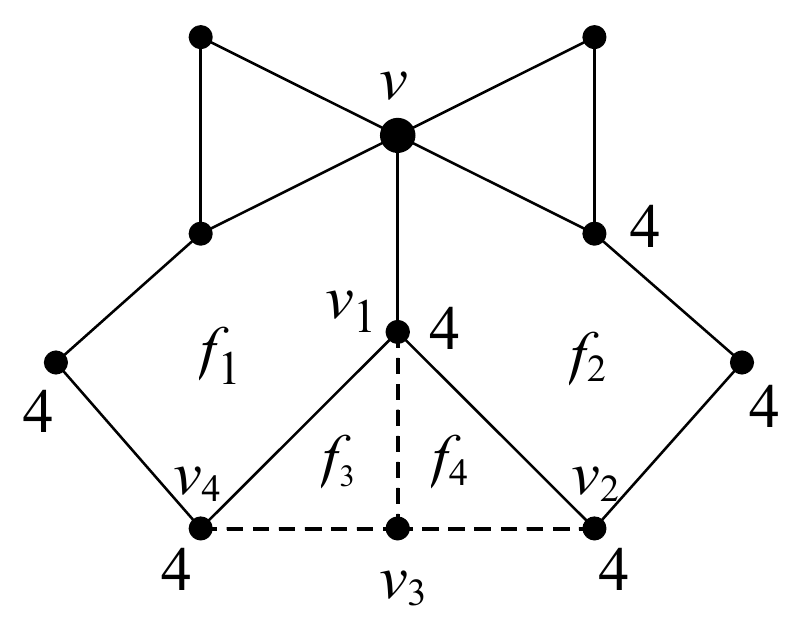}}
 \end{center}
  \caption{Illustration of Lemma \ref{bad5face1}} \label{fig2-4}
\end{figure}
\Proof
Assume that two bad $5$-faces, say $f_1$ and $f_2$, are normally adjacent, and they have one common $5^+$-vertex $v$ and one common edge $vv_1$ (see Figure \ref{fig2-4}). Since
both $f_1$ and $f_2$ are bad, then $v_1$ is a $4$-vertex, and $f_3=[v_1v_3v_4v_1]$ and $f_4=[v_1v_2v_3v_1]$ are two $3$-faces. Hence $f_3$ and $f_4$ are adjacent, contrary to (b). Thus the assumption is false.
\qed


\subsection{Discharging}

It follows that
\begin{align}
\sum_{v\in V}(2d_G(v)-6)+\sum_{f\in F}(d_G(f)-6)=-12 \nonumber
\end{align}
from Euler's formula $|V|-|E|+|F|=2$ and the equality $\sum_{v\in V}d_G(v)=2|E|=\sum_{f\in F}d_G(f)$.
Now we define an initial charge function $ch(x)$ for each $x\in V\cup F$ by letting $ch(v)=2d_G(v)-6$ for each $v\in V$ and  $ch(f)=d_G(f)-6$ for each $f\in F$. We are going to design several discharging rules. Since the sum of total charge is fixed during the discharging procedure, if we can change the initial charge function $ch(x)$ to the final charge function $ch'(x)$ such that $ch'(x)\ge0$ for each $x\in V\cup F$, then
\begin{align}
0\le \sum_{x\in V\cup F}ch(x)=\sum_{x\in V\cup F}ch'(x)=-12,\nonumber
\end{align}
 which is a contradiction. It means that no counterexample to Theorem \ref{main-thm} exists.  Thus Theorem \ref{main-thm} holds.

\begin{figure}[htbp]
 \begin{center}
\scalebox{0.5}[0.5]{\includegraphics {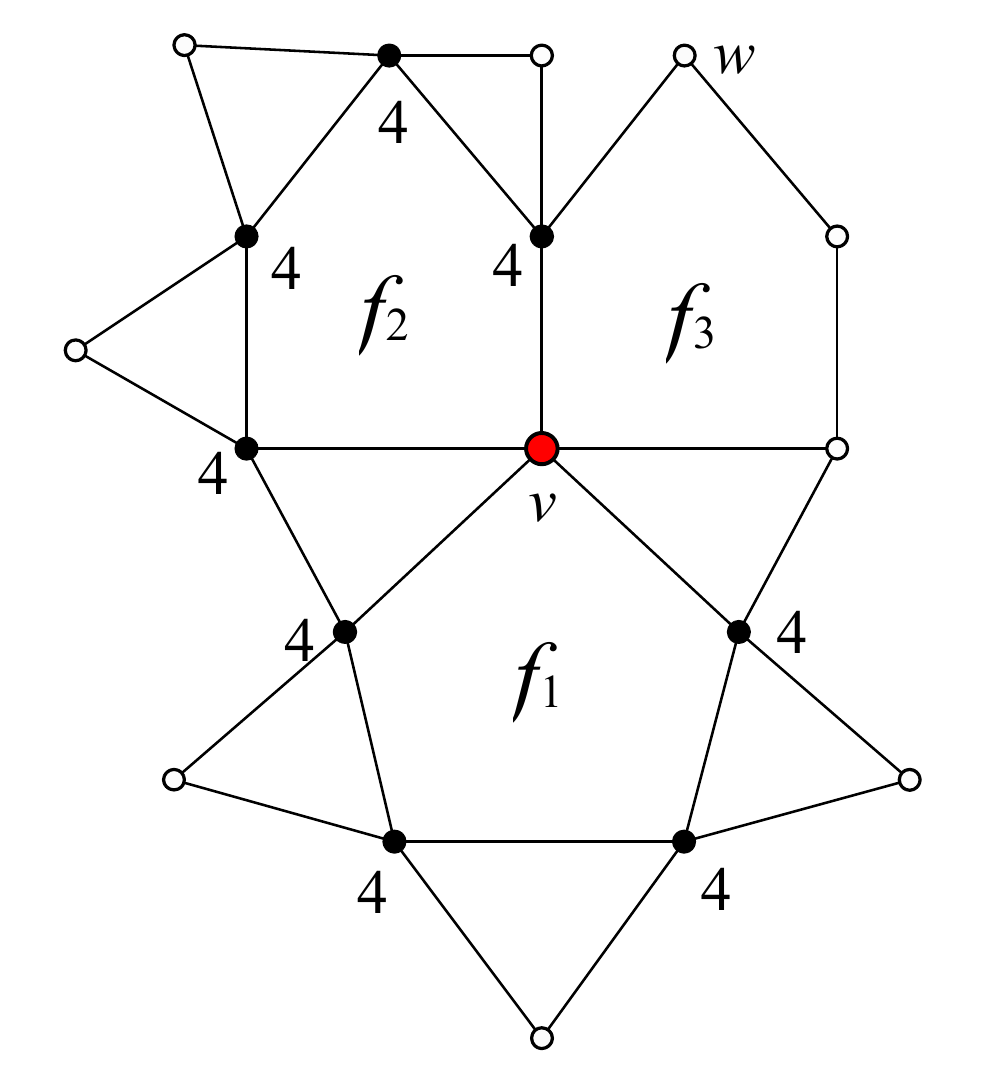}}
 \end{center}
 \caption{Special 5-vertex and incident 5-faces} \label{fig2-3}
\end{figure}

Before designing the discharging rules, we define a configuration as shown in Figure \ref{fig2-3}. In Figure \ref{fig2-3}, the degrees of white vertices are at least the number of edges incident to them, but the degrees of the black vertices are equal to the number of edges incident to them. The red vertex $v$ is a special $5$-vertex. In this configuration, $f_1$ is a special $5$-face, $f_2$ is a bad $5$-face, $f_3$ is a $5$-face satisfying the following:
\begin{itemize}
  \item $f_3$ is neither special nor bad; and
  \item $f_3$ is incident to a special $5$-vertex $v$, and $v$ is incident to a bad $5$-face $f_2$; and
  \item $f_3$ is normally adjacent to $f_2$.
\end{itemize}
We call $F_5$ the family of $f_3$, that is, $F_5$ is a family of $5$-faces that have the same properties of $f_3$.

\begin{remark}\label{remark3}
Let $w$ be a vertex on $f_3$ as shown in Figure \ref{fig2-3}. If $d_G(w)=4$, then $w$ is incident to at most one $3$-face according to (b).
\end{remark}

\bigskip


\noindent {\bf Discharging Rules:}

\begin{description}
\item[R1] Each $4^+$-vertex gives $1$ to each incident $3$-face.
\item[R2] Each $4^+$-vertex gives $\frac{1}{2}$ to each incident $4$-face.
\item[R3] Let $v$ be a $4$-vertex incident to at most one $3$-face. If $v$ is incident to exactly one $3$-face and one $4$-face, then $v$ gives $\frac{1}{4}$ to each incident $5$-face; Otherwise, $v$ gives $\frac{1}{3}$ to each incident $5$-face.
\item[R4] By Remark \ref{remark2} (1), a special $5$-face is exactly incident to one special vertex.
          \begin{description}
          \item[R4.1] Each special $5^+$-vertex gives $1$ to each incident special $5$-face.
          \end{description}
          \begin{description}

          \item[R4.2] Each special $5$-vertex gives $\frac{2}{3}$ to each incident bad $5$-face (see $f_2$ in Figure \ref{fig2-3}).
          \item[R4.3] If a special $5$-vertex is incident to three $5$-faces, which are special $5$-face $f_1$, bad $5$-face $f_2$, and non-special or non-bad $5$-face $f_3$ respectively, then $v$ gives $\frac{1}{3}$ to $f_3$. (See Figure \ref{fig2-3})
          \end{description}
\item[R5] Each non-special $5$-vertex and each $6^+$-vertex gives $\frac{3}{4}$ to each incident bad $5$-face.
\item[R6] Each source gives $\frac{1}{5}$ to each of its sinks.
\item[R7] If $u$ is a $6^+$-vertex, or $u$ is a non-special 5-vertex, or $u$ is a special $5$-vertex but $u$ is not incident to a bad $5$-face, then $u$ gives $\frac{1}{2}$ to each incident non-special or non-bad $5$-face.   
\end{description}

%
%
%

To complete the proof of Theorem \ref{main-thm}, it remains to check that the final charge of every element in $V\cup F$ is nonnegative. This will be shown by the following two claims.

\begin{claim}\label{claim1}
It holds that $ch'(v)\ge0$ for all $v\in V(G)$.
\end{claim}

\Proof
According to Lemma \ref{basic-lemma}, there exists no $3^-$-vertex in $G$, thus we need to consider the $4^+$-vertices in the following.\\

\Case {1:} When $d_{G}(v)=4$.

In this case, we have that $ch(v)=4\times 2-6=2$. According to Lemma \ref{good vertex}, $v$ has no any sink.
If $v$ is incident to two $3$-faces, then $ch'(v)=ch(v)-2\times 1=0$ by R1.

And, if $v$ is incident to exactly one $3$-face, then $v$ is incident to at most one $4$-face according to (b). Assume that $v$ is incident to one $4$-face, then it is incident to at most two $5$-faces. Thus $ch'(v)\ge ch(v)-1-\frac{1}{2}-\frac{1}{4}\times 2=0$ by R1, R2 and R3. Otherwise, $v$ is incident to at most three $5$-faces. Therefore, $ch'(v)\ge ch(v)-1-\frac{1}{4}\times 3>0$ by R1 and R3.

If $v$ is not incident to any $3$-face, then $ch'(v)\ge ch(v)-\frac{1}{2}\times 4=0$ by R2 and R3.\\

\Case {2:} When $d_{G}(v)=5$.

In this case, we have that $ch(v)=5\times 2-6=4$.
Assume that $v$ is special, that is, $v$ is incident to a special $5$-face by Remark \ref{remark2} (1). From Lemma \ref{3-faces} and the definition of special $5$-face, it holds that $v$ is incident to exactly two $3$-faces and one special $5$-face. According to Lemma \ref{specialgood}, $v$ has no  sink. If $v$ is incident to a bad $5$-face, then $v$ is incident to at most one bad $5$-face by Lemma \ref{bad5face1}. We have that $ch'(v)\ge ch(v)-1\times 2-1-\frac{2}{3}-\frac{1}{3}=0$ from R1 and R4. Otherwise $v$ is not incident to any bad $5$-face. It holds that $ch'(v)\ge ch(v)-1\times 2-1-\frac{1}{2}\times 2=0$ by R1, R4 and R7.

We next assume that $v$ is not special, that is, $v$ is not incident to a special $5$-face by Remark \ref{remark2} (1). The vertex $v$ is incident to at most two $3$-faces by Lemma \ref{3-faces} (1). According to Lemma \ref{bad5face1}, $v$ is incident to at most two bad $5$-faces.
\begin{itemize}
  \item If $v$ is incident to two $3$-faces, then $v$ is incident to at most one bad $5$-face from Lemma \ref{bad5face1}, the definition of bad $5$-face and (b), and is not incident to any $4$-face by (b). Thus we have that $ch'(v)\ge ch(v)-1\times 2-\frac{3}{4}-\frac{1}{2}\times 2>0$ from R1, R4 and R7.
  \item If $v$ is incident to exactly one $3$-face, then $v$ is incident to at most two bad $5$-faces by the definition of bad $5$-face. Thus the final charge $ch'(v)> ch(v)-1-\frac{1}{5}-\frac{3}{4}\times 2-\frac{1}{2}\times 2>0$ from R1, R4, R6 and R7.
  \item If $v$ is not incident to any $3$-face, then $v$ is not incident to bad $5$-face. It holds that $ch'(v)\ge ch(v)-\frac{1}{2}\times 5>0$ by R2 and R7.
\end{itemize}

%
%
%

\Case {3:} When $d_{G}(v)=6$.

Note that we have that  $ch(v)=2\times 6-6=6$. According to Lemma \ref{3-faces}, $v$ is incident to at most three $3$-faces and at most three special $5$-faces.

\begin{itemize}
  \item Assume that $v$ is incident to three $3$-faces. If $v$ is incident to three special $5$-faces, then $v$ has no any sink by Lemma \ref{specialgood}. Thus we have that $ch'(v)\ge ch(v)-3-3=0$ from R1, R4.1. Otherwise $v$ is incident to at most two special $5$-faces. Then $v$ has no any sink by Lemma \ref{specialgood}.
  Hence we have that
$ch'(v)\ge ch(v)-3 - 2-\frac{3}{4}>0$
from R1, R4.1 and R5.
  \item Assume that $v$ is incident to at most two $3$-faces, then $v$ is incident to at most one special $5$-faces by Lemma \ref{3-faces} (2).
      If $v$ is incident to one special $5$-faces, then $ch'(v)> ch(v)-2(1+\frac{1}{5})-1-\frac{3}{4}\times3>0$ from R1, R4.1, R5 and R6. Otherwise,
      by Lemma \ref{bad5face1}, $v$ is incident to at most two bad $5$-faces. Therefore, it holds that $ch'(v)> ch(v)-2(1+\frac{1}{5})-\frac{3}{4}\times2-\frac{1}{2}\times 2>0$ from R1, R4.1, R5, R6 and R7.
\end{itemize}

\Case {4:} When $d_{G}(v)= 7$.

In this case, we have that $ch(v)=2\times 7-6=8$. By Lemma \ref{3-faces}, $v$ is incident to at most three $3$-faces and at most two special $5$-faces. The smallest final charge $ch'(v)> ch(v)-3(1+\frac{1}{5})-2-\frac{3}{4}\times2>0$ from R1, R4.1, R5 and R6.\\

\Case {5:} When $d_{G}(v)= k\ge8$.

Observe that $ch(v)=2k-6$. By Lemma \ref{3-faces}, $v$ is incident to at most $\lfloor\frac{d_G(v)}{2}\rfloor$ $3$-faces and at most $\lfloor\frac{d_G(v)}{2}\rfloor$ special $5$-faces.
It holds that $ch'(v)\ge ch(v)-(1+\frac{1}{5})\times \lfloor\frac{d_G(v)}{2}\rfloor-1\times \lceil\frac{d_G(v)}{2}\rceil>0$ by R1, R4.1 and R6.
\qed

\begin{claim}\label{claim2}
It holds that $ch'(f)\ge0$ for all $f\in F(G)$.
\end{claim}

\Proof Let $f$ be a face of $G$. Because $G$ is simple, $G$ has no loops and multi-edges. Thus $d_G(f)\ge3$. If $d_G(f)\ge 6$, no charge is discharged from or to $f$, thus $ch'(f)=ch(f)=d_G(f)-6\ge0$.
If $d_G(f)=3$, then every vertex incident to $f$ gives $1$ to $f$ according to R1. Therefore, we have that $ch'(f)=ch(f)+3\times 1=d_G(f)-6+3=0$. If $d_G(f)=4$, then every vertex incident to $f$ gives $\frac{1}{2}$ to $f$ according to R2. Hence the final charge $ch'(f)=ch(f)+4\times \frac{1}{2}=d_G(f)-6+2=0$. Next we assume that $d_G(f)=5$. Note that $ch(f)=5-6=-1$.\\

\Case {1:} Assume that $f$ is small, that is, all the vertices incident to $f$ are $4$-vertices.

For $ 0 \leq t \leq 5$, let $t$ be the number of $4$-vertices which are incident to two 3-faces.  Then $f$ has $(5 - t)$ $4$-vertices which are incident to at most one 3-face, and $f$ has at least $t+1$ sources.  Thus we have
$
ch'(f)\ge -1+ \frac{1}{4}\times (5 - t) +\frac{1}{5}\times (t+1) = \frac{9 - t}{20} >0
$
for every $t \in \{0, 1, 2, 3, 4, 5 \}$ by R3 and R6.


\bigskip

\Case {2:} When $f$ is a $(5^+,4,4,4,4)$-face.

Denote the $5^+$-vertex by $v$. If $f$ is special, then $ch'(f)\ge -1+1=0$ according to R4.1. Next let $f$ be a non-special $5$-face.
\begin{itemize}
  \item Assume that $v$ is a special $5$-vertex and $v$ is incident to a bad
  $5$-face (see Figure \ref{fig2-3}). If $f$ is bad, then $ch'(f)\ge -1+\frac{2}{3}+\frac{1}{3}=0$ according to R4.2 and R3. Otherwise $f$ is not bad, then $f \in F_5$. Hence $f$ is incident to at least two $4$-vertices which are incident to at most one $3$-face, and are not incident to a 3-face and a 4-face at the same time by Remark \ref{remark3}, respectively. Therefore, we have that $ch'(f)\ge -1+\frac{1}{3}\times2+\frac{1}{3}=0$ according to R3, R4.3 and R7.

  \item Otherwise, $v$ is a $6^+$-vertex, or $v$ is not a special $5$-vertex, or $v$ is a special $5$-vertex and $v$ is not incident to a bad $5$-face. If $f$ is bad, then there exists a $4$-vertex incident to $f$, and the $4$-vertex is incident to at most one $3$-face. Thus we have that $ch'(f)\ge -1+\frac{1}{4}+\frac{3}{4}=0$ according to R3 and R5. Otherwise $f$ is not bad. Hence there are at least two $4^{+}$-vertices incident to $f$, and each of them is incident to at most one $3$-face. We can conclude that $ch'(f)\ge -1+\frac{1}{4}\times2+\frac{1}{2}=0$ according to R3 and R7.
\end{itemize}

\Case {3:} When there exist at least two $5^+$-vertices on $f$.
From Remark \ref{remark2} (2), we have that $f$ is neither special nor bad.
\begin{itemize}
  \item If $f \in  F_5$, then $ch'(f)\ge -1+\frac{1}{3}\times2+\frac{1}{3}=0$ by R3, R4.3 and R7.
  \item Otherwise $f\notin F_5$, that is, $f$ is not incident to a special $5$-vertex which is on a special $5$-face and is incident to a bad $5$-face. We can conclude that $ch'(f)\ge -1+\frac{1}{2}\times2=0$ according to R7. \qed
\end{itemize}

The proof of Theorem \ref{main-thm} is completed.
\qed

\section{Remarks}

\subsection{Difference between DP-coloring and list coloring}

A $\theta$-\emph{graph} is a graph consisting of two 3-vertices and three pairwise internally disjoint paths between the two 3-vertices. A $\theta$-\emph{subgraph} of $G$ is an induced subgraph that is isomorphic to a $\theta$-graph. We use $S\theta$ to denote such a \emph{special} $\theta$-\emph{subgraph} of $G$ in which one of the ends of the internal chord is a $5^-$-vertex and all of the other vertices are 4-vertices in $G$.

Let $G$ be a planar graph without 4-cycles adjacent to 3-cycles. In the proof of Theorem \ref{listversion} (A), the authors showed that if $G$ is not $4$-choosable with fewest vertices, then $G$ contains no subgraph isomorphic to $S\theta$ (see Lemma 4 in \cite{CCW}). But we can not claim that if $G$ is not DP-$4$-colorable with fewest vertices, then $G$ contains no subgraph isomorphic to $S\theta$. Next we give an explanation.

%

\begin{figure}[htbp]
 \begin{center}
\scalebox{0.65}[0.65]{\includegraphics {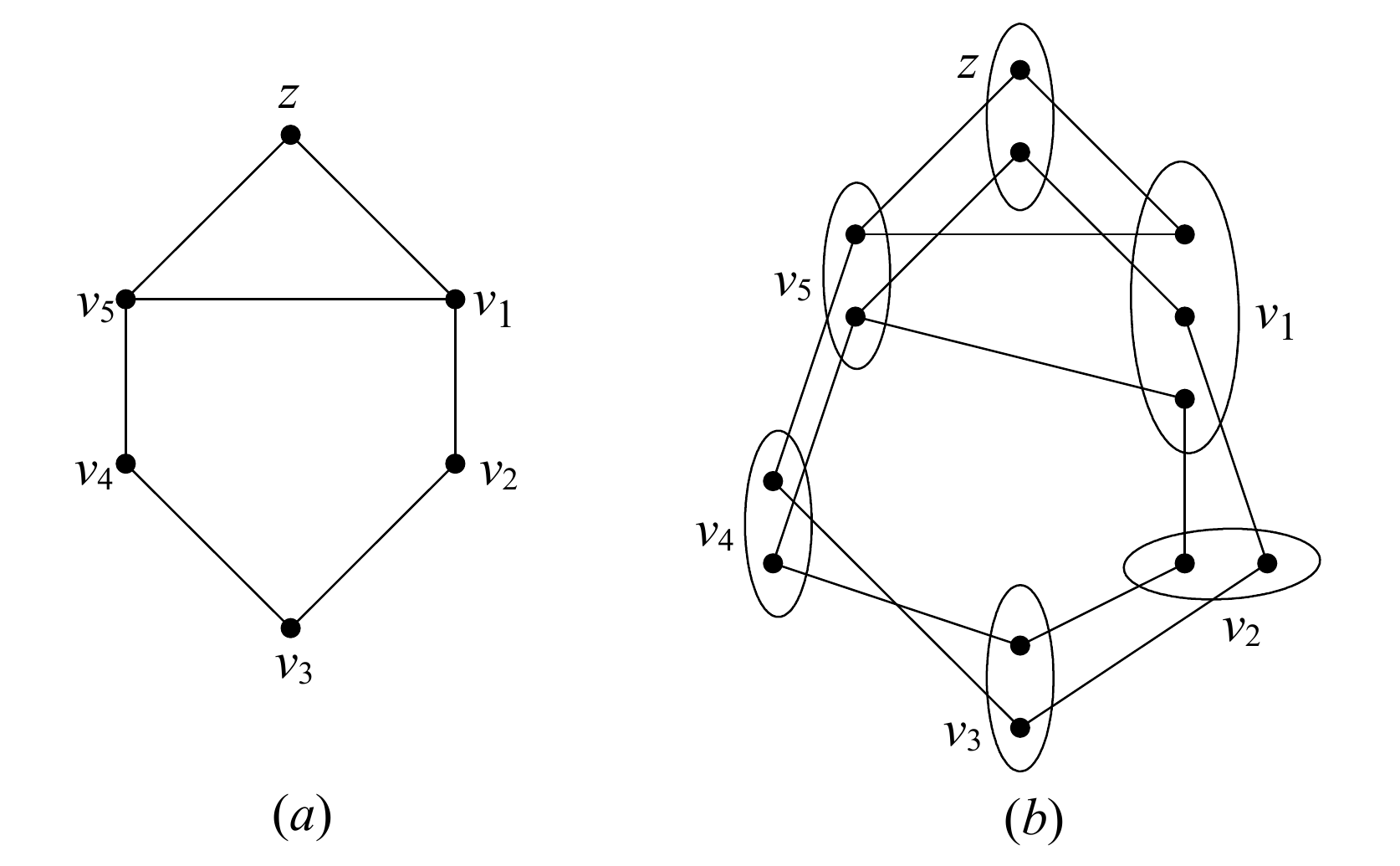}}
 \end{center}
 \caption{$S \theta$ graph} \label{S-theta}
\end{figure}

In Figure \ref{S-theta} (a), assume that $d_G(v_1) = 5$, $d_G(v_i) = 4$ for $i \in \{2, 3, 4, 5\}$ and $d_G(z) = 4$.
Let $G[S]$ be the subgraph of $G$ induced by $S = \{z, v_1, v_2, v_3, v_4, v_5 \}$.
Let $L$ be a list assignment of $G$ with $|L(v)|\ge 4$ for all $v\in V$,
and let $\M_{L}$ be a matching assignment over $L$.
Set $G' := G - S$ and $L'(v) = L(v)$ for $v \in V(G')$.
By (d), $G'$ admits an $\M_{L'}$-coloring.
Thus there is an independent set $I'$ in the $\M_{L'}$-cover with $|I'| = |V(G)| - |S| = |V(G)| - 6$.  For $v \in \{v_1, v_2, v_3, v_4, v_5, z \}$, we define
\[
L^*(v) = L(v) \setminus \bigcup_{uv \in E(G)} \big\{c' \in L(v): (u,c)(v,c') \in M_{L,uv}  \mbox{ and } (u, c) \in I' \big\}.
\]
Because $|L(v)| \ge 4 $ for all $v \in V(G)$,
we have that $|L^*(v_5)| \geq 3$, and $|L^*(v)| \geq 2$ for $v \in \{z,v_1,v_2, v_3, v_4\}$.
We denote by $\M_{L^*}$
the restriction of $\M_{L}$ into $G[S]$ and $L^*$.

Next we give a $\M_{L^*}$-cover as shown in (b) of Figure \ref{S-theta}. But we cannot find an independent set $I^*$ with $|I^*| = 6$ in $\M_{L^*}$-cover. Thus we cannot claim that if $G$ is not DP-$4$-colorable with fewest vertices, then $G$ contains no subgraph isomorphic to $S\theta$. Therefore, Theorem \ref{main-thm} is not trivial from Theorem \ref{listversion} (A).


\subsection{Relationship with signed coloring}
There is a concept of {\em signed coloring} of signed graphs, which
was first defined by Zaslavsky \cite{Zaslavsky}
with slightly different form,
and then modified by M\'{a}\v{c}ajov\'{a}, Raspaud, and \v{S}koviera \cite{MRS}
so that it would be a natural extension of an ordinary vertex
coloring.  For detail story about signed coloring, we refer readers to \cite{JKS, KO2, MRS}.

An interesting obervation in \cite{KO2} is that the signed coloring of a signed graph $(G,\sigma)$
is a special case of a DP-coloring of $G$.  Thus Theorem \ref{main-thm} implies the following corollary, which is an extension of the result in \cite{JKS}.

\begin{corollary}
A graph $G$ is signed $4$-choosable if $G$ is a planar graph without $4$-cycle adjacent to a $3$-cycle.
\end{corollary}



\begin{thebibliography}{99}



\bibitem{Alon} N.~Alon. Degrees and choice numbers.
Random Structures \& Algorithms 16 (2000), 364--368.

\bibitem{Bernshteyn} A.~Bernshteyn. The asymptotic behavior of the correspondence chromatic number.
Discrete Math. 339 (2016), 2680--2692.


\bibitem{BK2} A.~Bernshteyn, A.~Kostochka.
On differences between DP-coloring and list coloring.
arXiv:1705.04883, (2017), preprint.


\bibitem{BKP} A.~Bernshteyn, A.~Kostochka, S.~Pron.
On DP-coloring of graphs and multigraphs.
Sib. Math.l J. 58 (2017), 28--36.


\bibitem{BKZ} A.~Bernshteyn, A.~Kostochka, X.~Zhu.
DP-colorings of graphs with high chromatic number.
European J. Combin. 65 (2017), 122--129.

\bibitem{CCW} P. Cheng, M. Chen, Y. Wang. Planar graphs without 4-cycles adjacent to triangles are 4-choosable. Discrete Math. 339 (12) (2016), 3052--3057.

\bibitem{DP} Z.~Dvo\v{r}\'{a}k, L.~Postle.
Correspondence coloring and its application to list-coloring planar graphs without cycles of lengths 4 to 8.
J. Comb. Theory, Ser. B., (2017), In press.

\bibitem{ERT} P. Erdos, A.L.Rubin, H.Taylor. Choosability in graphs. Proc. West Coast Conf. on Combinatorics, Graph Theory and Computing, Congressus Numerantium XXVI (1979), 125-157.

\bibitem{JKS} L.~Jin, Y.~Kang, E.~Steffen.
Choosability in signed planar graphs.
Europ. J. Combin. 52 (2016), 234--243.


\bibitem{KO} S.-J.~Kim, K.~Ozeki.
A note on a Brooks' type theorem for DP-coloring.
arXiv:1709.09807, (2017), preprint.


\bibitem{KO2} S.-J.~Kim, K.~Ozeki.
A Sufficient condition for DP-4-colorability.
arXiv:1709.09809, (2017), preprint.


\bibitem{Lam}
P. C.-B. Lam, B. Xu, J. Liu.
The 4-Choosability of Plane Graphs without 4-Cycles. J. Comb. Theory, Ser. B 76(1) (1999), 117--126.

\bibitem{MRS} E.~M\'{a}\v{c}ajov\'{a}, A.~Raspaud, M.~\v{S}koviera.
The chromatic number of a signed graph.
the Electron. J Combin. 23 (2016), \#P1.14.



\bibitem{Thomassen} C. Thomassen. Every Planar Graph Is 5-Choosable. J. Comb. Theory, Ser. B. 62(1) (1994), 180--181.

\bibitem{Vizing} V.G.Vizing. Vertex colorings with given colors (in Russian). Diskret. Analiz., 29 (1976): 3-10.


\bibitem{Voigt} M.~Voigt.
A not 3-choosable planar graph without 3-cycles. Discrete Math. 146 (1995), 325--328.



\bibitem{Zaslavsky} T.~Zaslavsky.
Signed graph coloring.
Discrete Math. 39 (1982), 215--228.






\end{thebibliography}
\end{document}